\documentclass[reqno]{amsart}%[rm]{birkmult}%,graphicx,a4paper
\usepackage{amssymb}
\usepackage[pagebackref,hypertex]{hyperref}%,colorlinks=false
\allowdisplaybreaks[4]
\theoremstyle{plain}
\newtheorem{thm}{Theorem}
\newtheorem{lem}{Lemma}

\theoremstyle{remark}
\newtheorem{rem}{Remark}
\DeclareMathOperator{\td}{d}
\date{Completed on 24 November 2008 in VU's Student Village}
\date{Revised from Open-TJM-2003.tex on 28 November 2010}
\date{}

\begin{document}

\title[An inequality of the gamma and digamma functions]
{An inequality involving the gamma and digamma functions}

\author[F. Qi]{Feng Qi}
\address[F. Qi]{Department of Mathematics, College of Science, Tianjin Polytechnic University, Tianjin City, 300160, China; School of Mathematics and Informatics, Henan Polytechnic University, Jiaozuo City, Henan Province, 454010, China}
\email{\href{mailto: F. Qi <qifeng618@gmail.com>}{qifeng618@gmail.com}, \href{mailto: F. Qi <qifeng618@hotmail.com>}{qifeng618@hotmail.com}, \href{mailto: F. Qi <qifeng618@qq.com>}{qifeng618@qq.com}}
\urladdr{\url{http://qifeng618.wordpress.com}}

\author[B.-N. Guo]{Bai-Ni Guo}
\address[B.-N. Guo]{School of Mathematics and Informatics, Henan Polytechnic University, Jiaozuo City, Henan Province, 454010, China}
\email{\href{mailto: B.-N. Guo <bai.ni.guo@gmail.com>}{bai.ni.guo@gmail.com}, \href{mailto: B.-N. Guo <bai.ni.guo@hotmail.com>}{bai.ni.guo@hotmail.com}}

\begin{abstract}
In the paper, we establish an inequality involving the gamma and digamma functions and use it to prove the negativity and monotonicity of a function involving the gamma and digamma functions.
\end{abstract}



\thanks{The first author was supported partially by the Science Foundation of Tianjin Polytechnic University}

%\thanks{This paper was typeset using \AmS-\LaTeX}

\maketitle

\section{Introduction}

It is common knowledge that the classical Euler gamma function $\Gamma(x)$ may be defined for $x>0$ by
\begin{equation}\label{egamma}
\Gamma(x)=\int^\infty_0t^{x-1} e^{-t}\td t.
\end{equation}
The logarithmic derivative of $\Gamma(x)$, denoted by $\psi(x)=\frac{\Gamma'(x)}{\Gamma(x)}$, is called the psi or digamma function, and $\psi^{(k)}(x)$ for $k\in\mathbb{N}$ are called the polygamma functions. It is well known that these functions are fundamental and that they have much extensive applications in mathematical sciences.
\par
The aim of this paper is to establish an inequality involving the gamma and digamma functions. The result may be stated as the following theorem.

\begin{thm}\label{open-TJM-2003-lem2}
For $t\in(0,\infty)$, we have
\begin{equation}\label{gamma(t/(1+2t))}
\frac{1+2t}{2t^2}\biggl[\ln\Gamma\biggl(\frac{t}{1+2t}\biggr)-\ln\Gamma(t)\biggr]<1-\psi(t).
\end{equation}
\end{thm}

As applications of Theorem~\ref{open-TJM-2003-lem2}, the following negativity and monotonicity are obtained.

\begin{thm}\label{open-TJM-2003-thm3}
For $y\in\bigl(-1,-\frac12\bigr)$, the function
\begin{equation}\label{q(x,y)}
x\psi(x+y+1)-\ln\Gamma(x+y+1)+\ln\Gamma(y+1)-\frac{x^2}{2(y+1)(x+y+1)}
\end{equation}
is negative and decreasing with respect to $x\in\Bigl[-\frac{2(y+1)^2}{1+2y},\infty\Bigr)$.
\end{thm}

\begin{rem}
The negativity of the function~\eqref{q(x,y)} is equivalent to
\begin{equation}
\biggl[\frac{\Gamma(x+t)}{\Gamma(t)}\biggr]^{1/x}>\exp\biggl[\psi(x+t)-\frac{x}{2t(x+t)}\biggr]
\end{equation}
for $t\in\bigl(0,\frac12\bigr)$ and $x\in\bigl[-\frac{2t^2}{2t-1},\infty\bigr)$.
\end{rem}

\section{Lemmas}

In order to prove our main results, the following lemmas are needed.

\begin{lem}[{\cite[p.~305]{poly-comp-mon.tex}}]
For $x>0$, we have
\begin{equation}\label{qicui305}
\frac{1}{2x}-\frac{1}{12x^2}<\psi(x+1)-\ln x<\frac{1}{2x}.
\end{equation}
\end{lem}

\begin{lem}[{\cite[Lemma~1]{f-mean} and \cite[Theorem~1]{gamma-psi-batir.tex-jcam}}]
The inequality
\begin{equation}\label{correct-batir-ineq}
e^{\psi(L(a,b))}<\biggl[\frac{\Gamma(a)}{\Gamma(b)}\biggr]^{1/(a-b)}
\end{equation}
is valid for positive numbers $a$ and $b$ with $a\ne b$, where
\begin{equation}\label{L(a,b)}
L(a,b)=\frac{b-a}{\ln b-\ln a}
\end{equation}
stands for the logarithmic mean.
\end{lem}

\begin{lem}[{\cite[Theorem~1]{new-upper-kershaw-2.tex-mia}}]
For real numbers $s>0$ and $t>0$ with $s\ne t$ and an integer $i\ge0$, the inequality
\begin{equation}\label{new-upper-main-ineq}
(-1)^i\psi^{(i)}(L_p(s,t))\le \frac{(-1)^i}{t-s}\int_s^t\psi^{(i)}(u)\td u \le(-1)^i\psi^{(i)}(L_q(s,t))
\end{equation}
holds if $p\le-i-1$ and $q\ge-i$, where $L_p(a,b)$ is the generalized logarithmic mean of order $p\in\mathbb{R}$ for positive numbers $a$ and $b$ with $a\ne b$,
\begin{equation}
L_p(a,b)=
\begin{cases}
\left[\dfrac{b^{p+1}-a^{p+1}}{(p+1)(b-a)}\right]^{1/p},&p\ne-1,0;\\[1em]
\dfrac{b-a}{\ln b-\ln a},&p=-1;\\[1em]
\dfrac1e\left(\dfrac{b^b}{a^a}\right)^{1/(b-a)},&p=0.
\end{cases}
\end{equation}
\end{lem}

\begin{lem}\label{qi-psi-ineq-beta-lem}
For $x\in(0,\infty)$ and $k\in\mathbb{N}$, we have
\begin{equation}\label{qi-psi-ineq-beta-1}
\ln\biggr(x+\frac12\biggl)-\frac1x<\psi(x)<\ln(x+1)-\frac1x
\end{equation}
and
\begin{equation}\label{qi-psi-ineq-beta-2}
\frac{(k-1)!}{(x+1)^k}+\frac{k!}{x^{k+1}}< (-1)^{k+1}\psi^{(k)}(x) <\frac{(k-1)!}{(x+1/2)^k}+\frac{k!}{x^{k+1}}.
\end{equation}
\end{lem}

\begin{proof}
In \cite[Theorem~1]{Guo-Qi-Srivastava2007.tex}, the following necessary and sufficient conditions are obtained: For real numbers $\alpha\ne 0$ and $\beta$, the function
\begin{equation*}
g_{\alpha,\beta}(x)=\biggl[\frac{e^x\Gamma(x+1)} {(x+\beta)^{x+\beta}}\biggr]^\alpha,\quad x\in(\max\{0,-\beta\},\infty)
\end{equation*}
is logarithmically completely monotonic if and only if either $\alpha>0$ and $\beta\geq1$ or $\alpha<0$ and $\beta\leq\frac12$. Further considering the fact in~\cite[p.~98]{Dubourdieu} that a completely monotonic function which is non-identically zero cannot vanish at any point on $(0,\infty)$ gives
$$
(-1)^k[\ln g_{\alpha,\beta}(x)]^{(k)}=(-1)^k\alpha[x+\ln\Gamma(x)+\ln x-(x+\beta)\ln(x+\beta)]^{(k)}>0
$$
for $k\in\mathbb{N}$ and $x\in(0,\infty)$ if and only if either $\alpha>0$ and $\beta\geq1$ or $\alpha<0$ and $\beta\leq\frac12$. As a result, from straightforward calculation and standard arrangement, inequalities \eqref{qi-psi-ineq-beta-1} and~\eqref{qi-psi-ineq-beta-2} follow. Lemma~\ref{qi-psi-ineq-beta-lem} is thus proved.
\end{proof}

\begin{lem}[{\cite[p.~296]{kuang-3rd} and \cite[p.~274, (3.6.19)]{mit}}]
For $x>0$, we have
\begin{equation}\label{kuang-mit-ineq-ln}
\ln\biggl(1+\frac1x\biggr)<\frac2{2x+1}\biggl[1+\frac1{12x}-\frac1{12(x+1)}\biggr].
\end{equation}
\end{lem}

\begin{lem}[{\cite[p.~107, Lemma~3]{theta-new-proof.tex-BKMS}}]\label{comp-thm-1}
For $x\in(0,\infty)$ and $k\in\mathbb{N}$, we have
\begin{equation}\label{qi-psi-ineq-1}
\ln x-\frac1x<\psi(x)<\ln x-\frac1{2x}
\end{equation}
and
\begin{equation}\label{qi-psi-ineq}
\frac{(k-1)!}{x^k}+\frac{k!}{2x^{k+1}}< (-1)^{k+1}\psi^{(k)}(x)<\frac{(k-1)!}{x^k}+\frac{k!}{x^{k+1}}.
\end{equation}
\end{lem}

\begin{rem}
It is noted that the left-hand side inequality~\eqref{new-upper-main-ineq} for $i=0$ and $p=-1$ is just the inequality~\eqref{correct-batir-ineq}. For more information about inequalities~\eqref{correct-batir-ineq} and~\eqref{new-upper-main-ineq}, please refer to \cite{bounds-two-gammas.tex, new-upper-kershaw-JCAM.tex, subadditive-qi-guo-jcam.tex} and related references therein.
\end{rem}

\begin{rem}
In~\cite[Corollary~3]{egp} and \cite[Theorem~1]{Infinite-family-Digamma.tex}, it was proved that the double inequality
\begin{equation}\label{corollary2.3-rew}
    \ln\biggl(x+\frac12\biggr)-\frac1x< \psi(x)< \ln(x+e^{-\gamma})-\frac1x
\end{equation}
holds on $(0,\infty)$. In \cite[Theorem~1]{Infinite-family-Digamma.tex}, it was also shown that the scalars $\frac12$ and $e^{-\gamma}=0.56\dotsm$ in~\eqref{corollary2.3-rew} are the best possible. It is obvious that the inequality~\eqref{corollary2.3-rew} refines and sharpens~\eqref{qi-psi-ineq-beta-1}.
\end{rem}

\begin{rem}
In~\cite{Sharp-Ineq-Polygamma.tex}, the inequality~\eqref{qi-psi-ineq-beta-2} was refined and sharpened.
\end{rem}

\section{Proofs of theorems}

Now we are in a position to prove our theorems.

\begin{proof}[Proof of Theorem~\ref{open-TJM-2003-lem2}]
It is clear that
\begin{equation*}
\frac{1+2t}{2t^2}\biggl[\ln\Gamma\biggl(\frac{t}{1+2t}\biggr)-\ln\Gamma(t)\biggr] =-\frac{1+2t}{2t^2}\int_{t/(1+2t)}^t\psi(u)\td u,
\end{equation*}
so the required inequality~\eqref{gamma(t/(1+2t))} can be rewritten as
\begin{equation*}
p(t)\triangleq\frac{1+2t}{2t^2}\int_{t/(1+2t)}^t\psi(u)\td u-\psi(t)+1>0,\quad t>0.
\end{equation*}
The inequality~\eqref{qicui305} is equivalent to
\begin{equation}\label{psi(x+1)ineq-var}
  \ln x-\frac{1}{2x}-\frac{1}{12x^2}<\psi(x)<\ln x-\frac{1}{2x},\quad x>0.
\end{equation}
From~\eqref{psi(x+1)ineq-var}, it follows that
\begin{align*}
p(t)&>\frac{1+2t}{2t^2}\int_{t/(1+2t)}^t\biggl(\ln u-\frac{1}{2u}-\frac{1}{12u^2}\biggr)\td u
-\biggl(\ln t-\frac{1}{2t}\biggr)+1\\
&=\frac{4t-3\ln(2t+1)-1}{12t^2}\\
&\triangleq\frac{q(t)}{12t^2}.
\end{align*}
Since $q'(t)=\frac{2(4t-1)}{2t+1}$, the function $q(t)$ increases for $t>\frac14$. By $q\bigl(\frac87\bigr)=0.002\dotsm>0$, it is easy to see that the inequality~\eqref{gamma(t/(1+2t))} holds for $t\ge\frac87=1.1428\dotsm$.
\par
Letting $a=t$ and $b=\frac{t}{1+2t}$ in~\eqref{correct-batir-ineq} and~\eqref{L(a,b)} gives
\begin{equation*}
\frac{\ln\Gamma(t/(1+2t))-\ln\Gamma(t)}{t/(1+2t)-t}>\psi\biggl(\frac{2t^2}{(1+2t)\ln(1+2t)}\biggr).
\end{equation*}
Since the inequality~\eqref{gamma(t/(1+2t))} can be rearranged as
\begin{equation*}
\frac{\ln\Gamma(t/(1+2t))-\ln\Gamma(t)}{t/(1+2t)-t}>\psi(t)-1
\end{equation*}
for $t>0$, it is sufficient to show that
\begin{equation}\label{suffice-1}
\psi(t)-\psi\biggl(\frac{2t^2}{(1+2t)\ln(1+2t)}\biggr)
=\int_{2t^2/(1+2t)\ln(1+2t)}^t\psi'(u)\td u<1
\end{equation}
holds on $\bigl(0,\frac87\bigr)$.
\par
Taking $s=\frac{2t^2}{(1+2t)\ln(1+2t)}$, $i=1$ and $p=-2$ in the left-hand side of the inequality~\eqref{new-upper-main-ineq} leads to
$$
\int_{2t^2/(1+2t)\ln(1+2t)}^t\psi'(u)\td u \le \biggl[t-\frac{2t^2}{(2t+1)\ln(2t+1)}\biggr] \psi'\biggl(\biggl[\frac{2t^3}{(2t+1)\ln(2t+1)}\biggr]^{1/2}\biggr).
$$
Combining this with~\eqref{suffice-1} reveals that it suffices to prove
\begin{equation}\label{suffice-2}
\psi'\biggl(\biggl[\frac{2t^3}{(2t+1)\ln(2t+1)}\biggr]^{1/2}\biggr) \le\frac{(2t+1)\ln(2t+1)}{t[(2t+1)\ln(2t+1)-2t]}
\end{equation}
for $t\in\bigl(0,\frac87\bigr)$.
\par
The right-hand side inequality in~\eqref{qi-psi-ineq-beta-2} for $k=1$ results in
\begin{multline*}
\psi'\biggl(\biggl[\frac{2t^3}{(2t+1)\ln(2t+1)}\biggr]^{1/2}\biggr) \le\frac{(2t+1)\ln(2t+1)}{2t^3} \\
+\frac{1}{\sqrt{{2t^3}/{(2t+1)\ln(2t+1)}}\,+{1}/{2}}.
\end{multline*}
Then, in order to prove~\eqref{suffice-2}, it is enough to show
\begin{multline}\label{suffice-3}
\frac{(2t+1)\ln(2t+1)}{2t^3}+\frac{1}{\sqrt{{2t^3}/{(2t+1)\ln(2t+1)}}\,+{1}/{2}} \\
\le\frac{(2t+1)\ln(2t+1)}{t[(2t+1)\ln(2t+1)-2t]}
\end{multline}
for $t\in\bigl(0,\frac87\bigr)$.
\par
The inequality~\eqref{kuang-mit-ineq-ln} can be rewritten as
\begin{equation}
\ln(1+t)<\frac{t\bigl(t^2+12t+12\bigr)}{6(t+1)(t+2)},\quad t>0.
\end{equation}
Therefore, to verify~\eqref{suffice-3}, it is sufficient to prove
\begin{multline}\label{triangleq-q(t)-rew}
\frac{(2t+1)}{2t^3}\cdot\frac{2t\bigl(4t^2+24t+12\bigr)}{6(2t+1)(2t+2)} +\frac{1}{\sqrt{\frac{2t^3}{2t+1}\cdot\frac{6(2t+1)(2t+2)}{2t(4t^2+24t+12)}}\,+\frac{1}{2}}\\*
\le\frac1{t\Bigl[1-\frac{2t}{2t+1}\cdot\frac{6(2t+1)(2t+2)}{2t(4t^2+24t+12)}\Bigr]}
\end{multline}
for $t\in\bigl(0,\frac87\bigr)$, which can be simplified as
\begin{gather}
\frac{t^2+6t+3}{3t^2(t+1)}+\frac{2}{1+\sqrt{{12t^2(t+1)}/(t^2+6t+3)}\,}
\le\frac{t^2+6t+3}{t^2(t+3)},\notag\\
\sqrt{\frac{12t^2(t+1)}{t^2+6t+3}}\,\ge\frac{3 t^3+11 t^2+3 t-3}{t^2+6 t+3} \triangleq\frac{q(t)}{2\bigl(t^2+6t+3\bigr)}.\label{triangleq-q(t)}
\end{gather}
Since $q(t)$ is increasing on $(0,\infty)$ with $q(0)=-3$ and $q(1)=14$, the function $q(t)$ has a unique zero $t_0\in(0,1)$. From $q\bigl(\frac13\bigr)=-\frac23$, we can locate more accurately that $t_0\in\bigl(\frac13,1\bigr)$. When $0<t\le t_0$, the function $q(t)$ is non-positive, so the inequality~\eqref{triangleq-q(t)} is clearly valid. When $t\ge t_0$, the function $q(t)$ is non-negative, so squaring both sides of~\eqref{triangleq-q(t)} and simplifying gives
\begin{equation*}
h(t)\triangleq9t^6+54t^5+55t^4-60t^3-93t^2-18t+9\le0.
\end{equation*}
Direct differentiation yields
\begin{align*}
h'(t)&=54t^5+270t^4+220t^3-180t^2-186t-18,\\
h''(t)&=270t^4+1080t^3+660t^2-360t-186,\\
h^{(3)}(t)&=1080t^3+3240t^2+1320t-360.
\end{align*}
It is clear that the function $h^{(3)}(t)$ is increasing with $\lim_{t\to\infty}h^{(3)}(t)=\infty$ and $h^{(3)}(0)=-360$, so the function $h^{(3)}(t)$ has a unique zero which is the unique minimum point of the function $h''(t)$. Since $h''(0)=-186$ and $\lim_{t\to\infty}h''(t)=\infty$, the function $h''(t)$ has a unique zero which is the unique minimum point of the function $h'(t)$. From $h'(0)=-18$ and $\lim_{t\to\infty}h'(t)=\infty$, we conclude that the function $h'(t)$ has a unique zero which is the unique minimum point of the function $h(t)$ on $(0,\infty)$. Due to $h(0)=9$, $h\bigl(\frac13\bigr)=-\frac{700}{81}$, $h\bigl(\frac87\bigr)=-\frac{404759}{117649}$ and $\lim_{t\to\infty}h(t)=\infty$, it is not difficult to see that the function $h(t)<0$ on $\bigl(\frac13,\frac87\bigr)$. As a result, the inequalities~\eqref{triangleq-q(t)}, and so~\eqref{triangleq-q(t)-rew}, holds on $\bigl(0,\frac87\bigr)$. The proof of Theorem~\ref{open-TJM-2003-lem2} is complete.
\end{proof}

\begin{proof}[Proof of Theorem~\ref{open-TJM-2003-thm3}]
Denote the function~\eqref{q(x,y)} by $q(x,y)$. Differentiating and using the right-hand side inequality in~\eqref{qi-psi-ineq} yields
\begin{align*}
\frac{\partial q(x,y)}{\partial x}&=x\biggl[\psi'(x+y+1)-\frac{(x+2y+2)}{2(y+1)(x+y+1)^2}\biggr]\\
&<x\biggl[\frac1{x+y+1}+\frac1{(x+y+1)^2}-\frac{(x+2y+2)}{2(y+1)(x+y+1)^2}\biggr]\\
&=x\biggl[\frac{x(1+2y)+2(y+1)^2}{2 (y+1) (x+y+1)^2}\biggr],
\end{align*}
the function $\frac{\partial q(x,y)}{\partial x}$ is negative for $(x,y)\in\Bigl[-\frac{2(y+1)^2}{1+2y},\infty\Bigr)\times\bigl(-1,-\frac12\bigr)$, so the function $q(x,y)$ is decreasing with respect to $x\in\Bigl[-\frac{2(y+1)^2}{1+2y},\infty\Bigr)$ for $y\in\bigl(-1,-\frac12\bigr)$.
\par
Furthermore, from
\begin{multline*}
q\biggl(-\frac{2(y+1)^2}{1+2y},y\biggr)=\frac{2(y+1)^2}{2y+1} \biggl[1-\psi\biggl(-\frac{y+1}{2y+1}\biggr)\biggr] -\ln\Gamma \biggl(-\frac{y+1}{2y+1}\biggr)+\ln\Gamma(y+1)\\
=\frac{2(y+1)^2}{2y+1} \biggl[1-\psi\biggl(-\frac{y+1}{2y+1}\biggr) +\frac{\ln\Gamma(y+1) -\ln\Gamma(-(y+1)/(2y+1))} {y+1-[-(y+1)/(2y+1)]}\biggr]
\end{multline*}
and the inequality~\eqref{gamma(t/(1+2t))}, it follows that the function $q\Bigl(-\frac{2(y+1)^2}{1+2y},y\Bigr)$ is negative for $y\in\bigl(-1,-\frac12\bigr)$, and so the function $q(x,y)$ is negative for $x\in\Bigl[-\frac{2(y+1)^2}{1+2y},\infty\Bigr)$ and $y\in\bigl(-1,-\frac12\bigr)$. The proof of Theorem~\ref{open-TJM-2003-thm3} is complete.
\end{proof}

\begin{rem}
This paper is a part of the preprint~\cite{Open-TJM-2003.tex}.
\end{rem}

\end{document}